\documentclass{article}

\usepackage[margin=3.5cm]{geometry}

\usepackage[utf8]{inputenc}
\usepackage{hyperref}
\usepackage{amsmath,amsfonts}
\usepackage{amssymb}
\usepackage{natbib}
\usepackage{comment}
\usepackage[parfill]{parskip}
\usepackage{graphicx}
\usepackage{authblk}
\usepackage{xcolor}

\usepackage{amsthm}

\newtheorem{theorem}{Theorem}
\newtheorem{corollary}{Corollary}

\newtheorem{procedure}{Procedure}
\newtheorem{proposition}{Proposition}

\newtheorem{fact}{Fact}

\DeclareMathOperator*{\dif}{d}

\usepackage{todonotes}

\title{Carefree multiple testing with e-processes}
\author[1]{Yury Tavyrikov}
\author[2]{Jelle J.\ Goeman}
\author[3]{Rianne de Heide\footnote{Corr.\ author: \texttt{r.deheide@utwente.nl}.}}
\affil[1]{Vrije Universiteit Amsterdam}
\affil[2]{Leiden University Medical Centre}
\affil[3]{University of Twente}

\date{\today} 

\begin{document}

\maketitle

\begin{abstract}
    E-processes enable hypothesis testing with ongoing data collection while maintaining Type I error control. However, when testing multiple hypotheses simultaneously, current $e$-value based multiple testing methods such as e-BH are not invariant to the order in which data are gathered for the different $e$-processes. This can lead to undesirable situations, e.g., where a hypothesis rejected at time $t$ is no longer rejected at time $t+1$ after choosing to gather more data for one or more $e$-processes unrelated to that hypothesis. We argue that multiple testing methods should always work with suprema of $e$-processes. We provide an example to illustrate that e-BH does not control the FDR, at level $\alpha$ when applied to suprema of $e$-processes. From the same example we see that the FWER is not controlled with averaging, and also closed e-BH does not control the FDR. We show that adjusters can be used to ensure FDR-sup control with e-BH under arbitrary dependence. 
\end{abstract}

\section{Introduction}

The $e$-process is a recent innovative conceptual framework for hypothesis testing. An $e$-process is a nonnegative stochastic process $(E_t)_{t \geq 0}$, which starts at a value of $E_0=1$ and is continually updated over time as more information comes in. It does not grow in expectation if the null hypothesis $H_0$ is true, that is, $\mathbb{E}_{P}(E_\tau)\leq 1$ for every stopping time $\tau$ and for every $P \in H_0$. Therefore, we may reject $H_0$ when, at any time point $t$, we have $E_t \geq 1/\alpha$. Doing so will control Type I error, since, if $H_0$ is true, Ville's inequality (\cite{ramdasHypothesisTestingEvalues2024}, Fact~6.4) guarantees that $e$-processes have the property that
\begin{equation} \label{eq:anytime_valid}
P\left(\sup_{t \in \mathbb{N}} E_t \geq 1/\alpha\right) \leq \alpha.
\end{equation}

In contrast to $p$-values, $e$-processes can be repeatedly evaluated until, hopefully, rejection follows for some $t$. This means that a researcher can continue gathering data as long as they like. Notably, where the definition of $e$-processes involves a stopping time, the Type I error control property \eqref{eq:anytime_valid} involves the supremum of the process. This implies that data gathering is carefree: the researcher never has to regret gathering more data. For example, suppose that the researcher decided to gather $E_{t+1}$ before properly evaluating $E_t$. If $E_t \geq 1/\alpha$ but $E_{t+1} < 1/\alpha$, then they can still reject $H_0$ and appeal to  \eqref{eq:anytime_valid} for Type I error control, effectively undoing the gathering of $E_{t+1}$, as long as $\alpha$ was fixed a priori (or even if it is not: see \citet{koning2023post, grunwald2024beyond}).

If more than one hypothesis is of interest, multiple testing considerations apply. There is by now a large literature on $e$-based multiple testing. Most notably, the e-BH procedure \citep{wangFalseDiscoveryRate2022} generalizes the Benjamini-Hochberg procedure for controlling the False Discovery Rate (FDR). However, this literature has so far almost exclusively focused on multiple testing with $e$-values, which are single instances of $e$-processes, i.e., random variables $E$ for which $\mathbb{E}(E)\leq 1$ if $H_0$ is true. When applied to parallel $e$-processes $E_t^1,\ldots,E_t^K$, as we will show below, e-BH guarantees that
\begin{equation} \label{eq:supFDR}
\sup_{t \in \mathbb{N}} \text{FDR}\left(E_t^1,\ldots,E_t^K\right) \leq \alpha.
\end{equation}

Unfortunately, data gathering under the guarantee \eqref{eq:supFDR} is not carefree. Let us illustrate this with an example. For $K=2$ $e$-values, e-BH rejects any hypothesis for which the $e$-value exceeds $2/\alpha$, and both hypotheses if both $e$-values exceed $1/\alpha$. Now suppose that at time $t$ we have $1/\alpha \leq E_t^1 < 2/\alpha$ and $E_t^2 < 1/\alpha$. Gathering more data for both $e$-processes, however, the situation could reverse, and $E_{t+1}^1 < 1/\alpha$ and $1/\alpha \leq E_{t+1}^2 < 2/\alpha$. Now, the researcher regrets gathering more data for $E^1$, since \eqref{eq:supFDR} does not support any rejection now, but would have supported rejection of both null hypotheses if the researcher would have stopped $E^1$, taking $E^1_{t+1}=E^1_t$. The lack of the carefree property complicates experimental design because collecting more data for one or more hypotheses can cause previously rejected hypotheses---those for which more data has been gathered, or even others---to lose their rejected status. This means that researchers may wish to consider carefully, in real time, for which hypotheses to gather more data, adding data where it seems most promising. As a result, the process becomes volatile and harder to manage. Carefree procedures avoid this problem by ensuring that collecting more data never leads to losing discoveries.

For carefree data gathering, the researcher must always be allowed to revert any $e$-process, for which they regret gathering more data, back to an earlier time point. The optimal reversion is to the maximum so far, so a carefree procedure should guarantee 
\begin{equation} \label{eq:FDRsup}
\sup_{s \in \mathbb{N}} \text{FDR}(\max_{t \leq s} E_t^1,\ldots, \max_{t \leq s} E_t^K) = \text{FDR}(\sup_{t\in \mathbb{N}} E_t^1,\ldots, \sup_{t\in \mathbb{N}} E_t^K) \leq \alpha.
\end{equation}
We argue that \eqref{eq:FDRsup}, the \emph{FDR-sup}, is the proper FDR criterion for $e$-processes. In this paper, we investigate how to construct procedures controlling \eqref{eq:FDRsup}.

\section{FDR control with e-processes}

\subsection{The e-BH procedure}

Let \(H_1, \dots, H_K\) be \(K\) hypotheses, denote the set of their indices with \(\mathcal{K} = \{1, \dots, K\}\), and let \(P^* \) be the true data-generating probability measure. For each \(k \in \mathcal{K}\), \(H_k\) is called \emph{true} if \(P^* \in H_k\).
Let \(\mathcal{N} \subseteq \mathcal{K}\) be the set of indices of true null hypotheses and let $K_0$ be its cardinality.
For each \(k \in \mathcal{K}\), \(H_k\) is associated with an $e$-value \(E^k\).

When there is only one observation of the $e$-value for each hypothesis, it has been shown \citep{wangFalseDiscoveryRate2022} that the classical Benjamini-Hochberg (BH) procedure can be modified to accommodate $e$-values. The most powerful property of this \emph{e-BH} procedure, is that it provides FDR control for any dependence structure between the $e$-values, in contrast with the classical BH procedure, which requires the assumption of positive regression dependence on a subset (PRDS). Both the classical BH and also the Benjamini-Yekutieli (BY) procedures can be recovered as special cases of the e-BH procedure.

\subsection{e-BH with e-processes}
Formally, a non-negative process $E$ is called an $e$-process when $\mathbb{E} [E_\tau] \leq 1$ for any stopping time $\tau$, where the expectation is taken under any distribution in the null hypothesis.

We consider the setting in which all $e$-processes are adapted to a common filtration $(\mathcal{F}_t)_{t \geq 0}$, as this is required to guarantee joint validity under stopping. This setting corresponds to many practical situations, and has been discussed in detail by \citet{wang2025anytimevalidfdrcontrolstopped}.

For each \(k \in \mathcal{K}\), let \(H_k\) be associated with an e-process \( (E^k_t)_{t\geq 0} \), and let
\(\tau_1, \tau_2, \ldots, \tau_K\) be their arbitrary stopping times. Let \(S_{k} = E^k_{\tau_k}\) be the stopped $e$-processes.

\begin{procedure}[Naive e-BH with e-processes] 
Let \(S_{[k]}\) for \(k \in \mathcal{K}\) be the $k$-th order statistic of \(S_1, \ldots, S_K\), sorted from the largest to the smallest, such that \(S_{[1]}\) is the largest stopped $e$-process.
The stopped e-BH procedure \(\mathcal{G}(\alpha): [0, \infty]^K \to 2^{\mathcal{K}}\) then rejects hypotheses with the largest \(k_{se}^*\) stopped $e$-processes, where
\[
k_{se}^* = \max \Bigg\{k \in \mathcal{K}: S_{[k]} \ge \frac{1}{\alpha} \frac{K}{k} \Bigg\}.
\]
\end{procedure}

\begin{fact}\label{fact:e-BHwitheprocesses}
The stopped e-BH procedure applied to arbitrarily dependent e-processes, adapted to a common filtration), has FDR at most \(K_0\alpha / K\).
\end{fact}

\noindent Since a stopped $e$-process has expectation under the null bounded by one, i.e.\ \(\mathbb{E} S_{k} = \mathbb{E} E^k_{\tau_k} \leq 1\), the set of stopped $e$-processes \( \{ S_{k} \}_{k \in \mathcal{K}}\) is a set of $e$-variables. Thus, this procedure reduces to the standard e-BH procedure on the stopped $e$-processes \( \{ S_{k} \}_{k \in \mathcal{K}}\), and the FDR guarantee \citep{wangFalseDiscoveryRate2022} carries over.

\subsection{Related work} Online multiple testing has been considered both classically with $p$-values \citep{fischerOnlineClosurePrinciple2024a} and $e$-values \citep{fischerOnlineClosedTesting2024, fischerOnlineGeneralizationEBH2024, xuOnlineMultipleTesting2023}, but online is here to be understood in the dimension of the hypothesis space: time points represent the addition of an hypothesis. Until recent, the only paper that considered multiple testing with $e$-values where data can be added sequentially, was \citep{xuOnlineMultipleTesting2023}, where they introduce the \emph{doubly sequential} framework. For $p$-values there is work on interim-analyses and group-sequential multiple testing, e.g.\ \citep{sarkar2019} on FDR; but to the best of our knowledge there is no work on any-time valid methods. Between uploading the present paper to arxiv and the submission, two notable developments appeared. \citet{wang2025anytimevalidfdrcontrolstopped} provide a rigorous setup for using e-BH on stopped e-processes. They raise the important issue that one needs to be careful about the filtrations used for the stopped e-processes. In the present paper we circumvented the issue by assuming a common underlying filtration, but it would be good to investigate the interplay between their result and ours. Secondly, two papers appeared in the same week: \citep{goeman2025epartitioningprinciplefalsediscovery}  and \citep{xu2025bringingclosurefdrcontrol}; proposing a novel necessary and sufficient principle for FDR controlling methods based on e-values. The authors have teamed up and are now writing a joint paper from the two preprints. Their principle gives rise to a substantial improvement of e-BH, we call it here \emph{closed e-BH}, and we added a subsection to show that our central proposition holds for closed e-BH as well. Lastly, \citet{hartog2025familywiseerrorratecontrol} study closed testing for familywise error rate (FWER) control, both in the static and sequential setting.

\section{Multiple testing with the running maximum of an e-process}
We have argued that multiple-testing procedures for $e$-processes should possesses the carefree property that sequential data gathering will not lead to regrets, and that this requires the procedure to work with the running maximum of the $e$-process: \(M^{k}_t = \max_{1 \le s \le t}E^k_{s}\). 

It is straightforward to observe that in such a procedure, the set of rejected hypotheses forms a non-decreasing set. That is, if a hypothesis is rejected at a given point in time, it will remain rejected in all subsequent time points. It therefore automatically possesses the accept-to-reject property, as introduced in \citep{fischerOnlineGeneralizationEBH2024}: a multiple testing method that can change earlier non-rejections into a rejection at a later point in time, but not the other way around. As noted before, in that work the sequential aspect was in terms of adding hypotheses rather than data points. We argue however that also in the latter context, the accept-to-reject principle is a desirable feature of any multiple testing method.

\begin{procedure}[Running maximum e-BH procedure]
Let \(M^{[k]}_t\) for \(k \in \mathcal{K}\) be the $k$-th order statistic of \(M^1_t, \ldots, M^K_t\), sorted from the largest to the smallest so that \(M^{[1]}_t\) is the largest maximum of an $e$-process.
The stopped e-BH procedure \(\mathcal{G}(\alpha): [0, \infty]^K \to 2^{\mathcal{K}}\) then rejects hypotheses with the largest \(k_{m}^*\) maxima of $e$-processes, where
\[
k_{m}^* = \max \Bigg\{k \in \mathcal{K}: M^{[k]}_t \ge \frac{1}{\alpha} \frac{K}{k} \Bigg\}.
\]
\end{procedure}

Since the running maximum of an $e$-process is not an $e$-process itself (see \cite{choeCombiningEvidenceFiltrations2024} for more on this), we cannot directly deduce FDR control for this procedure. However, in the case that all $e$-processes are independent, it is straightforward to see that this procedure ensures FDR control. This is because, in essence, the inverse of the running maximum of the $e$-processes \(1/M_k\) is a valid $p$-process \citep{choeCombiningEvidenceFiltrations2024}. Thus, applying this procedure to the running maximum is equivalent to applying the classical BH procedure.

However, in the case of arbitrary dependence between $e$-processes, the e-BH procedure loses its FDR control, as asserted by the following proposition.

\begin{proposition}\label{prop:FDRcontrolWithMaxE}
The running maximum e-BH procedure applied to arbitrarily dependent $e$-processes does not provide FDR control at level $\alpha$.
\end{proposition}

\subsection{Proof of Proposition~\ref{prop:FDRcontrolWithMaxE}}

We present a counterexample that demonstrates how the e-BH procedure, when applied to the running maxima of $e$-processes, fails to control the FDR at the nominal level \(\alpha\). The FDR is evaluated numerically.

Consider two true null hypotheses with corresponding $e$-processes \(E_t^1\) and \(E_t^2\), defined as:
\[
E_t^1 = X_0^1 \prod_{s=1}^t e_s^1, \quad E_t^2 = X_0^2 \prod_{s=1}^t e_s^2,
\]
where the initial values \((X_0^1, X_0^2)\) are distributed as:
\[
(X_0^1, X_0^2) \sim 
\begin{cases} 
(0, 0), & \text{with probability } 1 - 2\alpha, \\ 
\left(\frac{1}{2\alpha}, \frac{1}{2\alpha}\right), & \text{with probability } 2\alpha,
\end{cases}
\]
and the increments \((e_s^1, e_s^2)\) are distributed as:
\[
\begin{cases} 
\left(\frac{1}{2}, \frac{1}{2}\right), & \text{with probability } \frac{1}{3}, \\ 
\left(2, \frac{1}{2}\right), & \text{with probability } \frac{1}{3}, \\ 
\left(\frac{1}{2}, 2\right), & \text{with probability } \frac{1}{3}.
\end{cases}
\]
The variables \(X_0^1\), \(X_0^2\), and all \(e_s^1, e_s^2\) are independent. 

It is straightforward to verify that the expected values of \(E_t^1\) and \(E_t^2\) are equal to 1 at all times, i.e., \(\mathbb{E}[E_t^1] = \mathbb{E}[E_t^2] = 1\), which ensures that both \(E_t^1\) and \(E_t^2\) are valid $e$-processes. However, the processes \(E_t^1\) and \(E_t^2\) are not independent because the random variables \(e_s^1\) and \(e_s^2\), as well as \(X_0^1\), \(X_0^2\), are dependent.

Since the null hypotheses are true, the FDR-sup is equal to the probability of rejecting at least one of them. This corresponds to the event:
\begin{align}
    M_t^1 \geq \frac{2}{\alpha}, \quad M_t^2 \geq \frac{2}{\alpha}, \quad \text{or both } M_t^1, M_t^2 \geq \frac{1}{\alpha}, \label{eq:rejectionEventsFDR}
\end{align}

where \(M^{k}_t = \max_{1 \le s \le t}E^k_{s}\) as above.
To confirm the validity of this example, we performed simulations of the described e-processes (the code is provided in the Appendix~\ref{app:code}).

The results of the simulations show that, in this scenario, the FDR of the e-BH procedure applied to the running maxima is approximately \(1.08\alpha\). In our simulations, we used \( M = 10^6 \) Monte Carlo iterations. Therefore, based on the standard error formula \( \sqrt{p(1 - p) / M} \), the estimated standard deviation of the error is approximately \( 0.001 \). This demonstrates that when using running maxima, additional corrections or normalizations are necessary to ensure proper FDR control.

\subsection{Closed e-BH} Recently and in the same week, both \citet{goeman2025epartitioningprinciplefalsediscovery} and \citet{xu2025bringingclosurefdrcontrol} proposed a necessary and sufficient principle for FDR controlling methods. The principle gives rise to an improvement over standard e-BH, which we will now call \emph{closed e-BH}. In Theorem~2 of \citet{goeman2025epartitioningprinciplefalsediscovery}, it is proven that closed e-BH uniformly improves upon standard e-BH, and in particular, that means that the set of rejected hypotheses by e-BH is always included in the sets rejected by closed e-BH. Therefore, the proof of Proposition~\ref{prop:FDRcontrolWithMaxE} directly extends to closed e-BH.

\begin{corollary}
The procedure of applying closed e-BH to the running maxima of arbitrarily dependent e-processes does not control the FDR at level $\alpha$.
\end{corollary}

\subsection{FWER}
Recently, \citet{wang2024only} proved that the only admissible way of merging arbitrary e-values is to use a weighted arithmetic average, which implies that for FWER control under arbitrary dependence, we should average the e-values; mixing with the trivial e-value of $1$ if desired. We see from the events for which one or both hypotheses are rejected in e-BH, stated Equation~\ref{eq:rejectionEventsFDR}, that averaging those e-values will also always lead to rejection. Since e-values are always positive, if either one of them is larger than $2/\alpha$, the average will exceed $1/\alpha$, and if both of them are larger than $1\alpha$ the average is trivially larger than $1/\alpha$ as well. Therefore, Proposition~\ref{prop:FDRcontrolWithMaxE} extends to FWER as well.

\begin{corollary}
The procedure of averaging the running maxima of arbitrarily dependent e-processes does not control the FWER at level $\alpha$.
\end{corollary}

\section{Adjusting running maxima}
Since Proposition~\ref{prop:FDRcontrolWithMaxE} shows that e-BH does not control the FDR-sup, one way forward is to use \emph{adjusters}, which can transform the running maximum process such that it becomes an $e$-process. Adjusters were introduced in the context of option pricing \citep{dawidProbabilityfreePricingAdjusted2011} and game theoretic probability \citep{dawidInsuringLossEvidence2011, shaferTestMartingalesBayes2011}, ensuring against loss of capital or evidence, and in \citep{choeCombiningEvidenceFiltrations2024} as means to combine evidence across filtrations. 

An \emph{admissible adjuster} \cite{choeCombiningEvidenceFiltrations2024} is an increasing function $\mathrm{A} : [1, \infty] \rightarrow [0, \infty]$ if and only if it is right-continuous, $\mathrm{A}(\infty) = \lim_{E \rightarrow \infty} \mathrm{A}(E) = \infty$, and \mbox{
$\int_1^\infty \mathrm{A}(E) /E^2 \dif E = 1.$}
Examples can be found in the aforementioned papers, where they are named \emph{lookback adjuster, capital calibrator, martingale calibrator}, or \emph{adjuster}. If $M_t^k = \sup_{s\leq t} E_s$ is the running maximum of an $e$-process for a family of null hypotheses $H_0$, then $(\mathrm{A}(M_t^k))_{t\geq 0}$ is also an $e$-process for $H_0$.
Since the adjusted process is non-decreasing, it will never result in regretting having collected more evidence.

\subsection{FDR-sup control with the adjusted running maximum}

\begin{theorem}\label{theorem:adjustedeBH} The e-BH procedure applied to the adjusted running maxima of e-processes controls the FDR-sup at level $K_0 \alpha / K$.
\end{theorem}

 One way to prove this is to follow the classical proof of \cite{wangFalseDiscoveryRate2022}, and is detailed in Appendix~\ref{app:proofs}.

\subsection{Example: e-BH with the adjusted running maximum}
\begin{figure}[h]
\centering
\includegraphics[width=15cm, height=8cm]{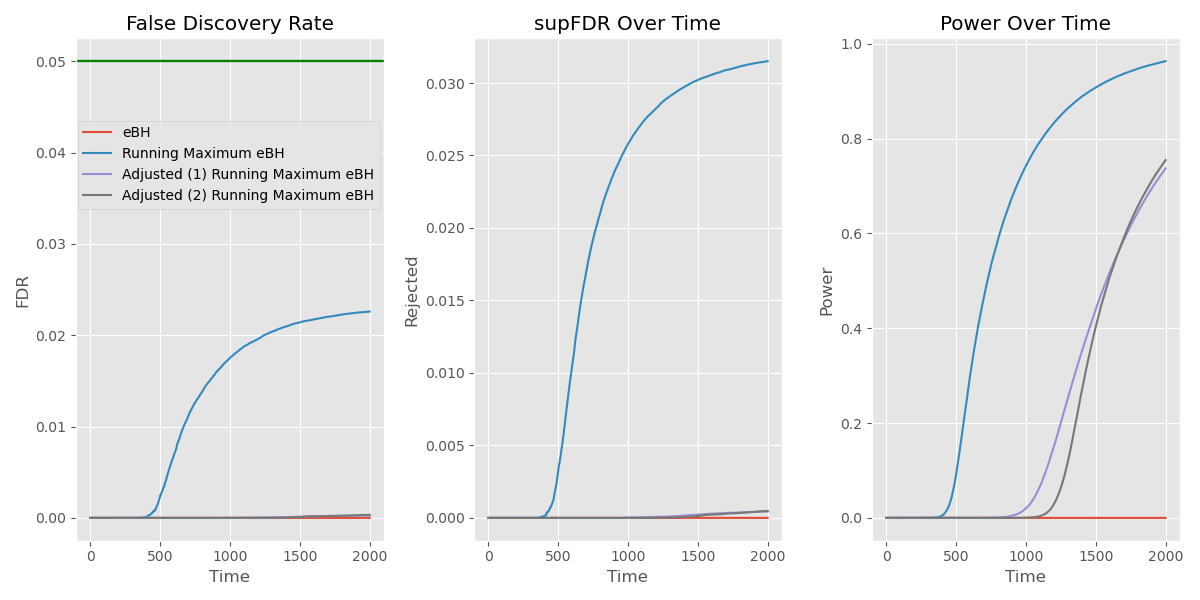}
\centering
\caption{Performance comparison of multiple testing methods (red line - the standard $e$-BH, blue line - running maximum $e$-BH, purple line - adjusted running maximum $e$-BH using $A_1$ from \eqref{eq:ExampleAdjuster}, black line - adjusted running maximum $e$-BH using $A_2$ from \eqref{eq:ExampleAdjuster}, green line - the confidence threshold). Plots: (1) FDR over time shows control below the threshold, (2) supFDR indicates worst-time FDR and (3) power - the proportion of rejected false hypotheses increases over time.}
\label{fig:sim}
\end{figure}

Two examples of admissible adjusters are
\begin{align}
\mathrm{A_1}(E) = \frac{E-1 - \log E}{\log^2E}, \quad\quad\quad \mathrm{A_2}(E) = \sqrt{E} -1.   \label{eq:ExampleAdjuster} 
\end{align}
In Figure~\ref{fig:sim} we see an example of running e-BH on the running maximum of an $e$-process adjusted by \eqref{eq:ExampleAdjuster}. The code to generate these plots can be found in Appendix~\ref{app:code}.

We conducted simulations to evaluate the performance of three procedures: standard e-BH, running maximum e-BH, and adjusted running maximum e-BH. The simulation involves $200$ hypotheses, where the null hypothesis (\(H_0\)) assumes a normal distribution with mean \(\mu = 0\), and the alternative hypothesis (\(H_1\)) assumes \(\mu = 0.1\). For each hypothesis, data is generated sequentially, creating a time series of observations.

The simulation includes dependencies among the hypotheses, introduced via a fixed randomly generated correlation matrix for the multivariate normal distributions for \(X^1_t, X^2_t, \ldots, X^K_t\), which are independent over time. To construct the $e$-process for each hypothesis, we compute the likelihood ratio at each time step as follows:
\[
E^k_t(X^t) = \prod_{i=1}^{t} \frac{f_1(X^k_i)}{f_0(X^k_i)} = \frac{e^{\frac{1}{2} \sum_{i=1}^t (X_i^k -0.1)^2}}{e^{\frac{1}{2}  \sum_{i=1}^t (X_i^k)^2}}.
\]

Each simulation generates \(T = 2000\) observations for each hypothesis, and the performance metrics (e.g., FDR, power, number of rejections) are averaged over \(M = 10000\) repetitions. This setup allows us to compare the efficiency and FDR control of the three methods under a dependent data structure.

In this example, we see that standard e-BH on the running maximum does not reach level $\alpha$. E-BH is known to be very conservative \citep{wangFalseDiscoveryRate2022, leeBoostingEBHConditional2024, xuMorePowerfulMultiple2024a}, which is why the methods seems to control FDR, even though Proposition~\ref{prop:FDRcontrolWithMaxE} says this is not the case in general. We also see that the adjusted maximum processes remain below level $\alpha$ as expected, and, unfortunately, we see an appreciable loss of power. Power could be boosted by using techniques such as conditional calibration \citep{leeBoostingEBHConditional2024} or stochastic rounding \citep{xuMorePowerfulMultiple2024a}. Still, the power loss of adjusted processes remains an important drawback of this way of controlling the FDR-sup. 


\section{Conclusion}

We have argued that multiple testing procedures based on 
$e$-processes should be carefree, meaning that the researcher should never lose out on any rejections because they chose to gather more data for one or more $e$-processes. If procedures are not carefree, experimental design becomes highly complex, since the decision when to gather more data for which $e$-processes becomes consequential. Carefree multiple testing procedures must operate on the running maxima of $e$-processes. Such running maxima are not themselves $e$-processes, but can be turned into e-processes by adjusters. Adjustment is costly, however, in terms of power. 


\paragraph{Acknowlegdgemenets} We thank Aaditya Ramdas and Nick Koning for helpful comments on an earlier version of this paper.

\paragraph{Declaration of funding} Rianne de Heide's work was supported by NWO Veni grant number VI.Veni.222.018. 

\bibliography{MaxEbib}

\appendix
\section{Proofs of FDR-sup control with adjusted e-BH.}
\label{app:proofs}

\subsection{Proof of Theorem~\ref{theorem:adjustedeBH}}
In one way to prove that the e-BH procedure with $e$-processes controls the FDR-sup, an adjuster is needed to make the last step of the classical proof by \cite{wangFalseDiscoveryRate2022} work. Let $\mathcal{R}^{\text{e-BH}}(\mathbf{e})$ denote the e-BH rejection set based on a set of $e$-values $\mathbf{e} = e^1_{t_1}, \ldots, e^K_{t_K}$, which can be obtained from different time points $t_1, \ldots, t_K$. I.e., $j \in \mathcal{R}^{\text{e-BH}}(\mathbf{e})$ if and only if $e^j_{t_j} > m / (\alpha \, | \mathcal{R}^{\text{e-BH}}(\mathbf{e}) |)$. Then the FDR-sup at time $s$ is
\begin{align}
    \text{FDR-sup}_s &= \sum_{j \in H_0} \frac{\mathbf{1} \left\lbrace \sup_{t\leq s} e^j_t \geq \frac{K}{\alpha |\mathcal{R}^{\text{e-BH}}(\mathbf{e})|} \right\rbrace }{|\mathcal{R}^{\text{e-BH}}(\mathbf{e})| \vee 1}\\ 
    &\leq \sum_{j \in H_0} \mathbb{E} \left[ \sup_{t \leq s} e^j_t \frac{\alpha\frac{ |\mathcal{R}^{\text{e-BH}}(\mathbf{e})|}{K}}{|\mathcal{R}^{\text{e-BH}}(\mathbf{e})| \vee 1} \right] \\ 
    &\leq \frac{\alpha}{K} \sum_{j \in H_0} \mathbb{E} \left[ \sup_{t \leq s} e^j_t \right].
\end{align}
The last expectation would be upper bounded by one, would the integrand be an $e$-value. This is however not the case, but when we use an adjuster $A$ on the expression, i.e.\ we take $A\left(\sup_{t \leq s}(e^j_t)\right)$, this \emph{is} an $e$-process, and the expectation is again upper bounded by one and the FDR-sup is controlled at level $ K_0\alpha / K \leq \alpha$.

\section{Code}\label{app:code}

\subsection{Code for the proof of Proposition~\ref{prop:FDRcontrolWithMaxE}}
The code can be found at the following link:  
\url{https://github.com/tavyrikov/runmax_eBH/blob/main/counterexpample.py}

\subsection{Code for generating Figure~\ref{fig:sim}}
The code can be found at the following link:  
\url{https://github.com/tavyrikov/runmax_eBH/blob/main/simulation.py}

\end{document}